\documentclass{article}
\usepackage{epsfig}
\usepackage{graphicx}
\usepackage[latin1]{inputenc}
\usepackage[T1]{fontenc}
\usepackage{mathptmx}
\usepackage{url}
\usepackage{amsmath, amssymb}
\newtheorem{definition}{Definition}
\newtheorem{theorem}[definition]{Theorem}
\newtheorem{lemma}[definition]{Lemma}
\newtheorem{remark}[definition]{Remark}
\newtheorem{corollary}[definition]{Corollary}

\newtheorem{example}{Example}

\begin{document}
\begin{center}
{\Large{On multivariate polynomials with many roots 
    over a finite grid\\
\ \\
}}
Olav Geil\\
Department of Mathematical Sciences\\
Aalborg University
\end{center}

\begin{abstract}
In this note we consider roots of multivariate polynomials
over a
finite grid. When given information on the leading monomial with
respect to a fixed monomial ordering, the footprint
bound~\cite{onorin, geilhoeholdt} provides us with an upper bound on the number
of roots, and this bound is sharp in that it can
always be attained by trivial polynomials being a constant times a product of an appropriate
combination of terms consisting of a variable minus a constant. In contrast
to the one variable case, there are multivariate polynomials
attaining the footprint bound being not of the above form. This even includes irreducible polynomials. The purpose of the note is to
determine a large class of polynomials for which only the mentioned
trivial polynomials can attain the bound, implying that to search for other
polynomials with the maximal number of roots one must look outside
this class.\\


\end{abstract}
\section{Introduction}

As is well-known, the number of roots of a univariate non-zero polynomial $F(X)$ over a
field ${\mathbb{F}}$ is upper bounded by its 
degree, and  if $a \in {\mathbb{F}}$ is a root then $X-a $ divides
$F$. For multivariate polynomials the situation is very different in
that such polynomials often possess infinitely many roots
when the field is not finite, and that $(a_1, \ldots , a_m)$ can be a
root of $F(X_1, \ldots , X_m)$ without $F$ being divisible by any
$X_\ell-a_\ell$, $\ell \in \{1, \ldots , m\}$. 

In many applications, however,
the point set under consideration is finite, e.g. it corresponds to a
grid $S_1 \times \cdots \times S_m$ where $S_\ell$,
 $\ell=1, \ldots , m$ are finite subsets of ${\mathbb{F}}$, one particular
case being the grid ${\mathbb{F}}_q^m$ where ${\mathbb{F}}_q$ denotes the finite
field with $q$ elements. For finite grids, of course, the number of
roots of a polynomial $F$ becomes finite and the Schwartz-Zippel lemma
provides us with an upper bound using information on the total degree of $F$.
Having knowledge of the leading monomial with respect to some fixed
monomial ordering, the footprint bound produces even more precise
information, as we recall in a moment. 

One is often 
interested in determining polynomials which attain the maximal number
of roots according to the above bounds. For instance this is the case
when one wants to produce good algebraic geometric
codes~\cite{stichtenoth,handbook} or when trying to determine minimal weight code words of generalized
Reed-Muller codes and their relatives,
e.g.~\cite{grmdef,massey1973polynomial,abs,hyperbolic,ruano2009structure,ruano2007parameters,weighted,hiram}. As it turns out a way to produce
polynomials with the maximal number of roots according to the
footprint bound is to take a non-zero constant times a product of
terms of the form
$X_\ell-a$, but as already hinted this trivial construction will typically not give us
all the desired polynomials. For instance it does not produce any
irreducible polynomials of degree more than one. In the present note
we determine a large
class of polynomials for which the above type of trivial polynomials are
exactly those attaining the maximal number of roots, implying that to
search for other polynomials with the maximal number of roots one must look
outside this class.

The note is organized as follows. In Section~\ref{secfoot} we recall
the footprint bound for a single polynomial and demonstrate its
sharpness. Then in Section~\ref{secmain} we give necessary conditions
for non-trivial polynomials to attain the footprint bound.

\section{The footprint bound}\label{secfoot}

The footprint bound~\cite{onorin,geilhoeholdt} is a general method to upper bound the size of
varieties~\cite{clo}. For a single polynomial $F(X_1, \ldots , X_m)$ and the
point set being a finite grid one obtains a simple closed formula
expression. This formula assumes that $\deg_{X_\ell} {\mbox{lm}}(F) < s_\ell$, $\ell=1, \ldots ,
m$, where $s_\ell$ denotes the size of $S_\ell$ and where
${\mbox{lm}}(F)$ denotes the leading monomial of $F$. Of course this
condition in particular is satisfied when not only the leading monomial, but all
monomials in the support satisfy the conditions on the degree which we
will often assume throughout the note. Observe
that even the last assumption is no real restriction as $F(X_1, \ldots ,
X_m)$ has exactly the same roots over the finite grid~\cite{alon} as does 
\begin{equation}
F(X_1, \ldots ,
X_m) {\mbox{ rem }} \big\{\prod_{\alpha \in S_1}(X_1-\alpha), \ldots ,
\prod_{\alpha \in S_m} (X_m-\alpha)\big\},\label{eqremainder}
\end{equation} 
where the latter notation means the
remainder modulo the polynomials in the curly brackets. This remainder is 
produced using the multivariate division
algorithm~\cite{clo}. It is clear that~(\ref{eqremainder}) satisfies the
requirement on the degree for  
each monomial in the support of it and we therefore introduce the
following notation for the set of remainders:
\begin{eqnarray}
 &&{\mbox{\  \hspace{-1cm} \ }} {\mathbb{F}}[X_1, \ldots , X_m]_{<(s_1, \ldots , s_m)}= \nonumber
  \\
&& \{ F(X_1,
                                                            \ldots ,
                                                            X_m) \in
                                                            {\mathbb{F}}[X_1,
                                                            \ldots ,
                                                            X_m] \mid  \deg_{X_\ell} F < s_\ell, {\mbox{\  for \ }} \ell =1, \ldots , m\}. \nonumber
\end{eqnarray}
The footprint bound for a single polynomial now is~\cite{onorin,geilhoeholdt}:
\begin{theorem}\label{thefoot}
Let the notation be as above and consider an arbitrary, but fixed,
monomial ordering on the set of monomials in the variables
$X_1, \ldots , X_m$. For a non-zero polynomial  $F(X_1, \ldots , X_m)
\in {\mathbb{F}}[X_1, \ldots , X_m]$ write its leading monomial as
$X_1^{i_1} \cdots X_m^{i_m}$ and assume $i_1 < s_1, \ldots , i_m
<s_m$ (e.g.\ $F(X_1, \ldots , X_m) \in {\mathbb{F}}[X_1, \ldots ,
X_m]_{<(s_1, \ldots , s_m)}$). Then $F$ possesses at most $D(X_1^{i_1} \cdots X_m^{i_m})=s_1 \cdots s_m - (s_1-i_1) \cdots
(s_m-i_m)$ roots over $S_1 \times \cdots \times S_m$. 
\end{theorem}
We remark that the Schwartz-Zippel
lemma~\cite{ore1922uber,demillo,schwartz,zippel} can be obtained as a corollary. This bound states that
over a finite grid $S\times \cdots \times S$, a polynomial of total
degree $t <s$ can at most have $ts^{m-1}$ roots where $s=\# S$.

 It is
not difficult to see that Theorem~\ref{thefoot} is sharp in that for
any prescribed leading monomial corresponding polynomials exist having as many
roots as the upper bound. Namely, consider any subsets
\begin{equation}
S_\ell^\prime \subseteq S_\ell, i_\ell =\# S_\ell^\prime, {\mbox{ \ }} \ell=1, \ldots ,
m \label{eq1} 
\end{equation}
and $k \in {\mathbb{F}} \backslash \{0\}$. Then the trivial polynomial 
\begin{equation}
k \prod_{\ell=1}^m\prod_{\alpha \in S_\ell^\prime}(X_\ell-\alpha) \label{eq2}
\end{equation}
has exactly $(s_1-i_1) \cdots (s_m-i_m)$ non-roots in the finite grid
and therefore the number of roots as predicted in
Theorem~\ref{thefoot}, as obviously for any monomial ordering the
leading monomial of (\ref{eq2}) equals $X_1^{i_1}\cdots
X_m^{i_m}$. 
 Observe that the polynomial in~(\ref{eq2}) is an
example of a polynomial satisfying that if $(a_1, \ldots , a_m) \in
S_1 \times \cdots \times S_m$ is a root then some $X_\ell-a_\ell$ divides
it. 

The polynomials described in~(\ref{eq2}) are by no means the only ones producing equality in
the footprint bound. We illustrate this observation with a classic example.
\begin{example}\label{ex1}
Let $q$ be a prime power and consider the point set
${\mathbb{F}}_{q^2} \times {\mathbb{F}}_{q^2}$. The Hermitian
polynomial $F(X_1, X_2)= X_1^{q+1}-X_2^q-X_2$ has $q^3$ roots which is exactly the
upper bound from Theorem~\ref{thefoot} when choosing a monomial
ordering such that $X_2^q$ becomes the leading monomial. The roots are
established by using the fact that $X_1^{q+1}$ is the norm function
related to the field extension 
${\mathbb{F}}_{q^2}/{\mathbb{F}}_q$ and that $X_2^q+X_2$ is the
similar trace function. Employing the properties of these functions one
determines~\cite{stichtenothhermitian,yanghermitian} the
roots. Clearly, for no $(a_1,a_2)\in {\mathbb{F}}_{q^2} \times {\mathbb{F}}_{q^2}$ it
holds that $X_1-a_1$, nor that $X_2-a_2$, divides
$F(X_1,X_2)$. Actually the Hermitian polynomial is absolutely irreducible.
\end{example}  

In the next section we derive information on when a polynomial can
possibly meet the footprint bound without being of the form~(\ref{eq2}).

\section{Necessary conditions for attaining the footprint bound}\label{secmain}

We start our investigations with a simple, yet crucial lemma.

\begin{lemma}\label{probiimp}
Let $F(X_1, \ldots , X_m) \in {\mathbb{F}}[X_1, \ldots , X_m]$ and $a
\in S_\ell$. Then the following
bi-implication holds true:
\begin{eqnarray}
&&X_\ell-a {\mbox{ \ divides \ }}  F(X_1, \ldots , X_m) \nonumber \\
& \Updownarrow \nonumber \\
&&(a_1, \ldots , a_{\ell-1}, a, a_{\ell+1}, \ldots , a_m) {\mbox{ is a root
   of }} F(X_1, \ldots , X_m) \nonumber \\
&&{\mbox{for all }} (a_1, \ldots , a_{\ell-1}, a_{\ell+1}, \ldots , a_m) \in
   S_1 \times \cdots \times S_{\ell-1} \times S_{\ell+1} \times \cdots
   \times S_m \nonumber
\end{eqnarray} 
\end{lemma}
{\bf{Proof:}} \ \ Without loss of generality we assume that $F(X_1,
\ldots , X_m) \in {\mathbb{F}}[X_1, \ldots , X_m]_{<(s_1, \ldots ,
  s_m)}$ (if this is not the case, we first perform reduction modulo $\{ \prod_{\alpha \in
  S_1}(X_1-\alpha), \ldots , \prod_{\alpha \in S_m}(X_m-\alpha)\}$). When $F$ is the zero polynomial then the bi-implication
clearly holds. Hence, assume that $F$ is not the zero polynomial. The "only if" part is easily verified. To to see the
"if" part we assume that $F$ has all the requested roots and apply the division algorithm to obtain $F(X_1, \ldots
, X_m)=Q(X_1, \ldots , X_{\ell-1},X_{\ell+1}, \ldots ,
X_m)(X_\ell-a)+R(X_1, \ldots , X_{\ell-1},X_{\ell+1}, \ldots ,
X_m)$. Aiming for a contradiction assume the $R$ is not the
zero-polynomial. Observe, that all monomials $M$ in the support of $R$ satisfy
$\deg_{X_i} M < s_i$, for all $i \in \{1, \ldots , \ell-1, \ell+1,
\ldots , m\}$ and in particular this holds for the leading monomial
with respect to any monomial ordering (the leading monomial exists by assumption). According to the footprint
bound applied to the point set $S_1\times \cdots \times
S_{\ell-1}\times S_{\ell+1} \times \cdots \times S_m$, $R$ therefore cannot have
all elements of this point set as roots. But $Q(X_1, \ldots ,
X_{\ell-1},X_{\ell+1}, \ldots , X_m)(X_\ell-a)$ has all the requested
roots of the lemma and therefore also this holds for $R$, which is a contradiction.   \hfill $\Box$ \\

We shall need two more lemmas.

\begin{lemma}\label{pro5}
Consider $F(X_1, \ldots , X_m) \in {\mathbb{F}}[X_1,
\ldots , X_m]$ with $\deg_{X_\ell} {\mbox{lm}}(F)<s_\ell$, $\ell=1,
\ldots , m$ with respect to some fixed monomial ordering (e.g.\ $F(X_1,
\ldots , X_m) \in {\mathbb{F}}[X_1, \ldots , X_m]_{<(s_1, \ldots ,
  s_m)}$). Write
$$F(X_1, \ldots , X_m)=G(X_1, \ldots , X_m)H(X_1, \ldots , X_m)$$ where
$$G(X_1, \ldots , X_m)=\prod_{i=1}^m \prod_{\alpha \in
  S_i^\prime}(X_i-\alpha)$$ for some $S_\ell^\prime \subseteq S_\ell$, $\ell=1,
\ldots , m$ and write $T_\ell=S_\ell \backslash S_\ell^\prime$ and $t_\ell=\# T_\ell$. Then $F$
attains the footprint bound over $S_1 \times \cdots \times S_m$ 
if and only if $H$ attains the footprint bound over $T_1 \times \cdots
\times T_m$.
\end{lemma}
{\bf{Proof:}} \ \ \\
The set of non-roots of $G$ over $S_1 \times \cdots \times S_m$ equals $T_1 \times \cdots \times
T_m$. Hence, the set of non-roots of $F$ over $S_1\times \cdots \times
S_m$ equals the set of  non-roots
of $H$ over $T_1 \times \cdots \times T_m$. Let
${\mbox{lm}}(F)=X_1^{i_1} \cdots X_m^{i_m}$, then
${\mbox{lm}}(H)=X_1^{i_1-s_1^\prime}\cdots X_m^{i_m-s_m^\prime}$ and 
according to the
footprint bound therefore the number of non-roots of $H$ over $T_1
\times \cdots \times T_m$ is at least
$$
\big( t_1-(i_1-s_1^\prime)\big) \cdots \big(
t_m-(i_m-s_m^\prime)\big)=(s_1-i_1)\cdots (s_m-i_m)$$ 
and exactly when equality holds the number of roots of $F$ over $S_1 \times \cdots
\times S_m$ becomes $s_1\cdots s_m-(s_1-i_1)\cdots (s_m-i_m)=D({\mbox{lm}}(F))$.
\ \hfill $\Box$ \\

\begin{example}\label{exflot}
Recall~\cite{reiter}, that the trace map ${\mbox{Tr}}: {\mathbb{F}}_{q^2} \rightarrow
{\mathbb{F}}_q$ is given by ${\mbox{Tr}}(\alpha)=\alpha^q+\alpha$ and
that the preimage of any element in ${\mathbb{F}}_q$ is of
size exactly $q$. Now let $F(X_1,X_2)=G(X_1,X_2)H(X_1,X_2) \in
{\mathbb{F}}_{q^2}[X_1,X_2]$ where 
$$G(X_1,X_2)=\prod_{{\mbox{Tr}}(\alpha)=0}(X_1-\alpha) \prod_{{\mbox{Tr}}(\alpha)=0}(X_2-\alpha)$$
and 
$$H(X_1,X_2)={\mbox{Tr}}(X_1)-{\mbox{Tr}}(X_2)=X_1^q-X_2^q+X_1-X_2.$$
We apply Lemma~\ref{pro5} to show that this
polynomial attains the footprint bound over $S_1 \times S_2$, where
$S_1=S_2={\mathbb{F}}_{q^2}$. The leading monomial of $F$ clearly is $X_1^{2q}X_2^q$ or $X_1^qX_2^{2q}$,
respectively, depending on the choice of monomial ordering. Hence,
according to the footprint bound $F$ potentially has
$q^4-(q^2-2q)(q^2-q)=3q^3-2q^2$ roots. The non-roots of $G(X_1,X_2)$ are 
$T_1 \times T_2$ where $T_1=T_2=\{\alpha \mid {\mbox{Tr}}(\alpha) \neq
0\}$. Therefore, according to Lemma~\ref{pro5}, $F$ attains the footprint bound over $S_1 \times S_2$ if and only if $H$ attains the footprint bound over
$T_1 \times T_2$ meaning that $H$ has
$(q^2-q)^2-(q^2-2q)(q^2-q)=q^3-q^2$ roots in
this set. This is exactly the case, the roots being 
$$\{(\alpha,\beta) \mid {\mbox{Tr}}(\alpha)={\mbox{Tr}}(\beta) \neq 0\}.$$
Hence, $F(X_1,X_2)$ does attain the footprint bound possessing
$3q^3-2q^2$ roots over $S_1 \times S_2={\mathbb{F}}_{q^2}\times {\mathbb{F}}_{q^2}$.
\end{example}

\begin{lemma}\label{lemlimlom}
Given the point set $S_1 \times \cdots \times S_m$ let $D$ be the
already introduced map and let $D^\prime$ be the  
map from $\{ X_1^{i_1} \cdots X_m^{i_m} \mid i_\ell <s_\ell, \ell =1,
\ldots , m\}$ to
${\mathbb{N}}_0$ defined by
$$D^\prime(X_1^{i_1}\cdots X_m^{i_m})=s_1\cdots s_{m-1}-(s_1-i_1) \cdots
(s_{m-1}-i_{m-1}).$$
Then for all monomials $M, N \in \{ X_1^{i_1} \cdots X_m^{i_m} \mid i_\ell <s_\ell, \ell =1,
\ldots , m\}$ we
have:
\begin{enumerate}
\item if $\deg_{X_m} M \geq 1$ then $D(M) > s_mD^\prime (M)$
\item if $N$ divides $M$ and $N\neq M$ then $D(N) < D(M)$. 
\end{enumerate}
\end{lemma}
{\bf{Proof:}} \ \ \\
By inspection. 
\ \hfill $\Box$ \\

We are now ready for our first result on which polynomials can
possibly attain the footprint bound.

\begin{theorem}\label{procond}\label{pro7}
Consider a non-constant polynomial  $H(X_1, \ldots , X_m)\in
{\mathbb{F}}[X_1, \ldots , X_m]_{<(s_1, \ldots , s_m)}$ with no factor
$X_\ell-a$, $a \in S_\ell$ for any $\ell$. Define $$\Omega = \max \{ D(M) \mid M \in
{\mbox{Supp}} (H) \}$$
where ${\mbox{Supp}}(H)$ denotes the support of $H$, and let $\{M_1, \ldots , M_\mu \}$ be the monomials in the support with $D$-value equal to $\Omega$. Then necessary conditions for $H$
to have $\Omega$ roots over $S_1\times \cdots \times S_m$ are that:
\begin{enumerate}
\item for any monomial ordering the leading monomial of $H$ belongs to
  $\{M_1, \ldots , M_\mu \}$
\item $\gcd (M_1, \ldots , M_\mu )=1$
\end{enumerate}
\end{theorem}
{\bf{Proof:}} \ \ \\
We only prove that condition 2 is needed. 
If $\gcd (M_1, \ldots , M_\mu ) \neq 1 $ then there exists some $X_\ell$
which divides $M_1, \ldots , M_\mu$. Without loss of generality 
assume $\ell=m$. For each $a \in S_m$ we then consider 
$$H^\prime(X_1,
\ldots , X_{m-1})=H(X_1, \ldots , X_{m-1},a)$$
which according to Lemma~\ref{probiimp} is non-zero by the
assumption that $H$ has no factor of the form $X_\ell-a$. 
Let $N$ be the leading monomial of $H^\prime$ with respect to
some monomial ordering on the set of monomials in the variable $X_1,
\ldots , X_{m-1}$. We now claim that $D^\prime (N) < \Omega /s_m$. If $N$ divides $M$ for some $M \in \{ M_1, \ldots ,
M_\mu \}$ then by 1.\ in Lemma~\ref{lemlimlom} we obtain $\Omega > s_m D^\prime
(M)$ which in turn is larger than or equal to
$s_m D^\prime (N)$ by 2.\ in the same lemma. If $N$ does not divide any monomial in $\{ M_1,
\ldots , M_\mu \}$ then for any $M$ in the support of $H$ such that
$N$ divides $M$ it holds that $s_m D^\prime (N) \leq s_m D^\prime
(M)\leq  D(M)$ which by assumption is
strictly smaller than $\Omega$. As $N$ necessarily divides some
monomial in the support of $H$ this proves the claim. By
Theorem~\ref{thefoot} for each $a$ the corresponding leading monomial
$N$ of $H^\prime$ gives us the upper bound $D^\prime(N)$ on the number of
roots of $H^\prime$ over $S_1 \times \cdots \times S_{m-1}$. Summing up the
contribution from each $a \in S_m$ we obtain fewer than
$\Omega$ roots in total of $H$ over $S_1 \times \cdots \times S_m$.    
\ \hfill $\Box$ \\

\begin{example}
This is a continuation of Example~\ref{exflot} where we considered the
polynomial $H(X_1,X_2)=X_1^q-X_2^q+X_1-X_2$ attaining the footprint
bound over $T_1 \times T_2$, with $T_1=T_2=\{ \alpha \in
{\mathbb{F}}_{q^2} \mid {\mbox{Tr}}(\alpha) \neq 0\}$. Clearly, $H$ has no
factor $X_\ell-a$, and therefore by 2.\ in  Theorem~\ref{procond} we
need to have
$\gcd(X_1^q,X_2^q)=1$ which indeed is the case. Furthermore, in
accordance with 1.\ in Theorem~\ref{procond} we have
$\Omega=D(X_1^q)=D(X_2^q)$. 
\end{example}

\begin{remark}\label{remrom}
If the polynomial $H$ in Theorem~\ref{procond} contains exactly
one monomial with the highest $D$-value then condition 2.\ is
never satisfied (here, we used the assumption from
Theorem~\ref{procond} that $H$ is not a constant). 
\end{remark}

\begin{theorem}\label{thetop}
Let $F(X_1, \ldots , X_m) \in {\mathbb{F}}[X_1, \ldots , X_m]_{<(s_1, \ldots , s_m)}$ be a
polynomial having a unique monomial $M$ in its support of highest
$D$-value and assume that for some monomial ordering this is the
leading monomial. Then it is either of the form~(\ref{eq2}) or it possesses fewer than $D(M)$ roots.
\end{theorem}
{\bf{Proof:}} \ \ \\
Assume that $F$ has $D(M)$ roots. Then $F$ is square-free and we may write $F=GH$ as in
Lemma~\ref{pro5}. We will show that $H$ is a constant. The crucial
observation is that
${\mbox{lm}}(H)$ is the unique monomial in the support of $H$ such
that the $D$-value with respect to the point set $T_1 \times \cdots
\times T_m$ is maximal. But then the assumption that $F$ has $D(M)$
roots by Lemma~\ref{pro5}, 2.\ in Theorem~\ref{pro7}, and
Remark~\ref{remrom} implies that $H$ is a constant.
\ \hfill $\Box$ \\

We present two corollaries to Theorem~\ref{thetop}. The first concerns
a family of polynomials which we call {\textit{monomial ordering invariant}}.

\begin{definition}\label{definvariant}
A polynomial is said to be monomial ordering invariant if it in its
support has a monomial $M$ which is divisible by any (other) monomial in
the support.
\end{definition}

Clearly, for a monomial ordering invariant polynomial $F$ and any choice
of monomial ordering the leading
monomial of $F$ equals the $M$ in
Definition~\ref{definvariant} and by 2.\ in Lemma~\ref{lemlimlom} the
$D$-value of $M$ is strictly larger than the $D$-value of any other monomial in the
support of $F$. As the name indicates there are no other
polynomials besides the monomial ordering invariant ones having a unique leading monomial. Let namely $M=X_1^{i_1}\cdots X_m^{i_m}$ be the
leading monomial of $F$ with respect to some monomial ordering and aiming for a contradiction assume that $F$
in its support has a polynomial $N=X_1^{i^\prime_1}\cdots
X_m^{i^\prime_m}$ with $i_t < i^\prime_t$ for some $t \in \{1, \ldots
, m\}$. But then for a lexicographic ordering with $X_t$ larger than the
other variables it holds that $N$ is larger than $M$ and
therefore $M$ is not the leading monomial with respect to this
ordering. Observe that univariate non-zero polynomials as well as the polynomials
in~(\ref{eq2}) satisfy the condition for being  monomial
ordering invariant. As an immediate corollary to Theorem~\ref{thetop}
we obtain:
\begin{corollary}\label{corolla}
Let $F(X_1, \ldots , X_m) \in {\mathbb{F}}[X_1, \ldots , X_m]_{< (s_1,
  \ldots , s_m)}$ be a monomial ordering invariant polynomial. Then
over $S_1 \times \cdots \times S_m$ the footprint bound is the same
for all choices of monomial ordering and when attained, $F(X_1,
\ldots , X_m)$ is necessarily of the form~(\ref{eq2}). 
\end{corollary}

\begin{remark}
Recall from the above discussion that a non-zero polynomial which is not
monomial ordering invariant must have in its support at least two
different monomials each being the leading monomial with respect to
some monomial ordering. But then according to the footprint bound if 
$F(X_1, \ldots , X_m) \in {\mathbb{F}}[X_1, \ldots , X_m]_{<(s_1, \ldots ,
  s_m)}$ is not monomial ordering invariant then the number of
roots cannot attain the maximal value of $D(M)$ for $M$ in the support
of it when the polynomial has only one monomial in its support with this $D$-value. Hence, the non-trivial information given in Theorem~\ref{thetop} is
actually that of Corollary~\ref{corolla}.
\end{remark}

The last corollary to Theorem~\ref{thetop} concerns irreducible
polynomials.

\begin{corollary}\label{cor2}
Let $F(X_1, \ldots , X_m) \in {\mathbb{F}}_q[X_1, \ldots ,
X_m]_{<(q,\ldots , q)}$ be an irreducible polynomial different from
$k(X_\ell-a)$, where $k \in {\mathbb{F}}_q \backslash \{0\}$, and with $s$ roots
from the point set $S_1 \times \cdots \times
S_m={\mathbb{F}}_q^m$. Then at least one of the following two
conditions are satisfied:
\begin{enumerate}
\item $F$ contains in its support at least two monomials with
  $D$-value equal to $s$
\item $F$ contains in its support a monomial of $D$-value strictly
  larger than $s$.
\end{enumerate}
\end{corollary}

\begin{example}
This is a continuation of Example~\ref{ex1} where we considered the
Hermitian polynomial $X_1^{q+1}-X_2^q-X_2$ which attains the footprint
bound, i.e.\ it possesses $D(X_2^q)=q^3$ roots over
${\mathbb{F}}_{q^2}$ (here we use a monomial ordering with $X_2^q$
larger than $X_1^{q+1}$). As already noted the Hermitian polynomial
is (absolutely) irreducible, and indeed it satisfies condition 2.\ of
Corollary~\ref{cor2} as $D(X_1^{q+1})=q^3+q^2>q^3$. 
\end{example}

\begin{remark}
One can interpret Theorem~\ref{procond},  Theorem~\ref{thetop} and Corollary~\ref{corolla} 
 as results on which polynomials different
from~(\ref{eq2}) can possibly attain the footprint bound over a finite grid.
\end{remark}

\section*{Acknowledgments}
The author would like to thank Ren\'{e} B{\o}dker Christensen for
comments on a previous version of the paper.

\end{document}